\documentclass[10 pt, leqno]{amsart}

\usepackage {amsfonts}
\usepackage{amsthm}
\usepackage{amssymb}
\usepackage{latexsym}
\usepackage{amsmath}
\usepackage{mathrsfs}

\theoremstyle{plain}
\newtheorem{tw}{Theorem}[section]

\newtheorem {prop}[tw] {Proposition}

\theoremstyle{definition}
\newtheorem {deft}[tw] {Definition}

\newcommand{\bn}{\Bbb N}

\newcommand{\alg} {\mathcal{A}}

\newcommand{\Hil}{\mathsf{H}}

\newenvironment{rlist}
{

\begin{enumerate}}
{\end{enumerate}}

\numberwithin{equation}{section}

\begin{document}

\author{Adam Skalski}
\address{Department of Mathematics and Statistics,  Lancaster University,
Lancaster, LA1 4YF} \email{a.skalski@lancaster.ac.uk}
\author{Joachim Zacharias}
\footnote{\emph{Permanent address of the first named author:} Department of Mathematics, University of \L\'{o}d\'{z}, ul. Banacha
 22, 90-238 \L\'{o}d\'{z}, Poland.}
\address{School of Mathematical Sciences,  University of Nottingham,
Nottingham, NG7 2RD}

 \email{joachim.zacharias@nottingham.ac.uk  }

\title{\bf A note on spectral triples and quasidiagonality}
\begin{abstract}
Spectral triples (of compact type) are constructed on arbitrary separable quasidiagonal $C^*$-algebras. On the other hand an
example of a spectral triple on a non-quasidiagonal algebra is presented.

\end{abstract}
 \maketitle

The concept of a spectral triple (unbounded Fredholm module) due to A.\,Connes  (\cite{book}) is a natural noncommutative
generalisation of a notion of a compact manifold, with certain summability properties corresponding in the classical case to the
dimension of the manifold. Recently E.\,Christensen and C.\,Ivan established the existence of spectral triples on arbitrary AF
algebras (\cite{ci}). In this note we generalise their result to arbitrary quasidiagonal (representations of) $C^*$-algebras.
Contrary to the AF situation our triples might in general have bad summability properties and it is not clear whether they
satisfy Rieffel's condition (i.e.\ whether the topology they induce on the state space coincides with the weak$^*$-topology).
Although the connections between properties related to quasi-diagonality and the existence of unbounded Fredholm modules seem to
have been known for a long time (see for example \cite{Voold}), explicit constructions have been until now given only in presence
of a filtration of the $C^*$-algebra in question consisting of finite-dimensional subspaces (\cite{Voold}, \cite{ci}).

 We also show that the existence of spectral triples of compact type does not imply quasidiagonality by exhibiting a
simple example of such a triple (with bad summability properties) on the natural, non-quasidiagonal, representation of the
Toeplitz algebra.

%It is known that the existence of spectral triples on $C^*$-algebras impose restrictions on the algebra in question.
%It is known that the existence
%of a $p$-summable spectral triple of traces and regularity conditions
%for the algebra question.

\section{Basics on quasidiagonality and spectral triples}

Throughout $A$ denotes a separable unital $C^*$-algebra. Representations of $A$ are assumed to act on separable Hilbert spaces.
% We write $\bn_0=\bn \cup \{0\}$.

\begin{deft}{(\cite{book})} A spectral triple or unbounded Fredholm module $(\alg, \Hil, D)$ on $A$ consists of a faithful representation
$\pi: A \to B(\Hil)$ together with a dense $^*$-subalgebra $\alg \subset A$ and an unbounded self-adjoint operator $D$ on $\Hil$
such that
\begin{rlist}
\item $[D,a]$ is densely defined and extends to a bounded operator on $\Hil$ for all $a \in \alg$;
\item $(I+D^2)^{-1}$ is compact.
\end{rlist}
If $p>0$ then $(\alg, \Hil, D)$ is $p$-summable if $(I+D^2)^{-p/2}$ is trace class. (Other summability conditions require
$(I+D^2)^{-1/2}$ to lie in various trace ideals.) Finally the triple $(\alg, \Hil, D)$ is \emph{of compact type} if $[D,a]$
defines a compact operator for all $a \in \alg$.
\end{deft}

It is known that the existence of spectral triples on $C^*$-algebras imposes restrictions on the algebra in question. For
instance the existence of a $p$-summable triple implies that $A$ is nuclear and has a tracial state (\cite{co}). Notice that the
Dirac operator on a compact spin  manifold $M$ of dimension $d$ defines a spectral triple on $C(M)$ which is $p$-summable for all
$p>d$ but not $p\leq d$ (see for example Chapter 7 in \cite{Roe}) .

\begin{deft}{(\cite{bo})}
A $C^*$-algebra $A$ is said to be quasidiagonal if there exists a sequence of completely positive and contractive maps
$\varphi_n:A \to M_{k_n}$ such that $\| \varphi_n(ab)- \varphi_n(a) \varphi_n(b) \| \to 0$ and $\| \varphi_n(a)\| \to \|a\|$ as
$n \to \infty$ for all $a,b \in A$.
\end{deft}

Every abelian $C^*$-algebra is quasidiagonal (use point evaluations on the spectrum); it is also easy to see that $AF$ algebras
are quasidiagonal. A representation $\pi: A \to B(\Hil)$ of a $C^*$-algebra is said to be quasidiagonal if there exists an
increasing sequence of finite rank projections $(P_n)_{n=1}^{\infty}$  on $\Hil$ such that $P_n$ converges strongly to $I$ and
$[\pi(a), P_n] \to 0$ as $n \to \infty$ for every fixed $a \in A$. D.\,Voiculescu showed in \cite{Voic} that a (separable)
$C^*$-algebra $A$ is quasidiagonal if and only if it admits a faithful quasidiagonal representation. Note that quasidiagonality
also implies the existence of tracial states (\cite{bo}, Proposition 7.1.16), if only $A$ is unital.

%Given a representation $\pi: A \to B(\Hil)$ of any $C^*$-algebra $A$

A quasidiagonal representation gives rise to the following setting. Given an increasing sequence $(P_n)_{n=1}^{\infty}$ of
projections on $\Hil$ converging strongly to $I$ for each $k \in \bn$ let $Q_k=P_k-P_{k-1}$, ($P_{0}:=0$). For $a \in B(\Hil)$
define the operator matrix $(a_{ij})_{i,j=1}^{\infty}$ by $a_{ij}=Q_i a Q_j$. Let moreover $(\alpha_i)_{i=1}^{\infty}$ be a
sequence of real numbers. Then we can define an essentially self-adjoint operator $D$ by  $D=\sum_{i=1}^{\infty} \alpha_i Q_i$
(with the sum understood strongly - we simply fix all eigenspaces and corresponding eigenvalues of $D$). If the projections $P_n$
have finite rank and $ |\alpha_i| \to \infty$ then $(I+D^2)^{-1}$ is compact. As $D$ is a diagonal operator with $\alpha_n$'s on
the diagonal it is clear that the matrix of $[D,a]$ is given by $((\alpha_i-\alpha_j)a_{ij})_{i,j=1}^{\infty}$:
$$
[D,a]=
\begin{pmatrix}
\alpha_1 a_{11} &\alpha_1 a_{12}  & \alpha_1 a_{13} & \dots & \\
\alpha_2 a_{21} &\alpha_2 a_{22}  & \alpha_2 a_{23} & \dots &\\
\alpha_3 a_{31} &\alpha_3 a_{32}  & \alpha_3 a_{33} & \dots &\\
\ldots &&&& \\
&&&& \\
\end{pmatrix}
-
\begin{pmatrix}
\alpha_1 a_{11} &\alpha_2 a_{12}  & \alpha_3 a_{13} & \dots & \\
\alpha_1 a_{21} &\alpha_2 a_{22}  & \alpha_3 a_{23} & \dots &\\
\alpha_1 a_{31} &\alpha_2 a_{32}  & \alpha_3 a_{33} & \dots &\\
\ldots &&&& \\
&&&& \\
\end{pmatrix}.
$$
In particular, if $[D,a]$ defines a bounded operator of norm $C$ then for all $i \neq j$ we must have $\|a_{ij}\| \leq
C|\alpha_i-\alpha_j|^{-1}$.

%As is well known it is difficult to

\section{Spectral triples on quasidiagonal $C^*$-algebras}

The following theorem implies in particular the existence of spectral triples on AF algebras, as proved in \cite{ci}. Contrary to
the AF situation  we cannot expect in this generality any good summability properties.

\begin{tw}
Let $A$ be a (separable) quasidiagonal $C^*$-algebra with quasidiagonal faithful representation $\pi: A \to B(\Hil)$ and let
$(b_i)_{i=1}^{\infty}$ be any sequence in $A$. Then there exists a spectral triple of compact type $(\alg, \Hil, D)$ on $A$, with
$\alg$ containing all $b_i$.
\end{tw}
\begin{proof}
By mixing our sequence $(b_i)_{i=1}^{\infty}$ with a dense sequence we obtain a dense sequence $(a_i)_{i=1}^{\infty}$; by taking
adjoints and finite products and putting them all in one sequence we may assume that  $(a_i)$ is moreover closed under taking
adjoints and products. This implies that $\text{span}(\{a_i : i \in \bn \})$ is a dense $^*$-subalgebra of $A$.

Now let $\pi: A \to B(\Hil)$ be a faithful quasidiagonal representation with a sequence of finite rank  projections
$(P_n)_{n=1}^{\infty}$ as before. Let $(\alpha_i)_{i=1}^{\infty}$ be an arbitrary sequence of real numbers such that $|\alpha_i|
\nearrow \infty$. Then writing $a$ for $\pi(a)$, where $a \in A$ we have (the sum is understood formally)
$$
\|[D,a]\| = \|\sum_k \alpha_k [Q_k,a] \|  \leq \sum_k |\alpha_k|( \|[P_k, a]\|+\|[P_{k-1},a]\|)
$$
so that $[D, a]= \sum_k \alpha_k [Q_k,a]$ converges in norm to a compact operator provided the right hand side converges. (Note
that $ [Q_k,a]$ is finite rank for all $k$.) All that remains to prove is the following statement:

\medskip

\noindent {\bf Claim:} There exists a subsequence of $(P_n)_{n=1}^{\infty}$ such that  $$\sum_k |\alpha_k|
(\|[a_i,P_{n_k}]\|+\|[a_i,P_{n_{k-1}}]\|)< \infty$$ for all $i \in \bn$.

\medskip

\noindent
\emph{Proof of Claim:}
Since $\|[a,P_k]\| \to 0$ for all $a \in A$ we can chose a subsequence $(P_{1,k})$ of $(P_n)$ such that
$$  \sum_k |\alpha_k| (\|[a_1,P_{1,k}]\|+\|[a_1,P_{1,k-1}]\|)< \infty.$$
Now chose a subsequence  $(P_{2,k})$ of $(P_{1,k})$ such that
$$  \sum_k |\alpha_k| (\|[a_2,P_{2,k}]\|+\|[a_2,P_{2,k-1}]\|)< \infty.$$
Since $|\alpha_k| \nearrow \infty$ we have
$$  \sum_k |\alpha_k| (\|[a_1,P_{2,k}]\|+\|[a_1,P_{2,k-1}]\|) \leq
 \sum_k |\alpha_k| (\|[a_1,P_{1,k}]\|+\|[a_1,P_{1,k-1}]\|)<\infty.$$
By induction we find a sequence of successive subsequences $(P_{l,k})$ such that
$$  \sum_k |\alpha_k| (\|[a_i,P_{l,k}]\|+\|[a_i,P_{l,k-1}]\|)< \infty$$
for $1\leq i \leq l$ and it is easy to see that the diagonal sequence $(P_{k,k})$ provides a required subsequence.
\end{proof}

Note that although we can choose the sequence $(\alpha_i)_{i=1}^{\infty}$ in an arbitrary way (as long as $|\alpha_i| \nearrow \infty$), the inductive construction above may entail that each $\alpha_i$ has very fast growing multiplicity in the list of eigenvalues of the operator $D$. This means
that unless we know some strong estimates on the rate of vanishing of the off-diagonal elements of elements of $A$ with respect
to the decomposition given by the original sequence $(P_n)_{n=1}^{\infty}$ we cannot expect the resulting triple to have any good
summability properties. When $A$ is an AF algebra then for any given $a \in A$ the off-diagonal elements with respect to the
natural sequence $(P_n)_{n=1}^{\infty}$ are simply $0$ from some point on, which explains why the triples constructed in
\cite{ci} can be arbitrarily well summable. For similar reasons we also cannot expect that the metric on the state space
$\mathcal{S}(A)$ given by the spectral triples constructed above induces the weak$^*$-topology on $\mathcal{S}(A)$ (in the spirit
of Rieffel's theory of compact quantum metric spaces, \cite{Rieffel}).

One might expect that the existence of a spectral triple $(\alg, \Hil, D)$ of compact type on a $C^*$-algebra $A$ should imply
that the representation $\pi:A \to B(\Hil)$ is quasidiagonal. This, however, is not true, as the next example shows:

\begin{prop}
Let $\pi:\mathcal{T} \to B(\Hil)$ be the standard (faithful) representation of the Toeplitz algebra $\mathcal{T}$ and let $\alg$
denote the $^*$-subalgebra of $\mathcal{T}$ generated by the unilateral shift. As $\mathcal{T}$ is not quasidiagonal, also the
representation $\pi$ is not quasidiagonal, but there exists a (non-finitely summable) spectral triple $(\alg, \Hil, D)$ of
compact type on $\mathcal{T}$.
\end{prop}
\begin{proof}
Consider the shift representation $\mathcal{T} = C^*(s) \to B(\ell^2(\bn))$, where $s$ is the unilateral shift and let
$\mathcal{P}(s)$ be the the $^*$-algebra generated by $s$. Then it is easy to check that $\mathcal{P}(s)=\text{span}\{1,
s^i,(s^*)^j, e_{i,j} : i, j \in \bn \}$ where $e_{i,j}$ denote the standard matrix units in $B(\ell^2(\bn))$. Let $(\alpha_i)$ be
a sequence of positive real numbers such that $\alpha_i  \nearrow \infty$ and $\alpha_{i+1}- \alpha_i \to 0$ as $i \to \infty$.
Then $[D,s]$ has the matrix representation
$$
\begin{pmatrix}
0 & 0 & \dots & & & \\
\alpha_{2}- \alpha_1 & 0 & \dots & && \\
0 &\alpha_{3}- \alpha_2 & 0 & \dots &&\\
0 &0 &\alpha_{4}- \alpha_3 & 0  &&\\
\ldots &&&&& \\
&&&&& \\
&&&&&
\end{pmatrix}
$$
which gives a compact operator and it is easy to see that $[D,a]$ is compact for every $a \in \mathcal{P}(s)$. However, $s$ and
$\mathcal{T}$ are clearly not quasidiagonal since $s$ has Fredholm index $-1$.

The triple constructed above will generally not be finitely summable, as can be seen for instance by putting $\alpha_n=
\sum_{k=1}^n k^{-1}$.
\end{proof}

We would like to finish the note with one more comment. If a $C^*$-algebra $A$ is residually finite dimensional (i.e.\ faithfully
embeddable into a  direct product of the form $\Pi_{i\in I} M_{k_i}$, see \cite{bo}) then one can trivially construct spectral triples of
arbitrarily good summability properties and such that
\begin{equation}\label{triv} [D,a]=0, \;\;\; a \in A.\end{equation} Such triples induce
the discrete topology on the state space of $A$ and thus do not characterise in the correct sense `topological dimension' of $A$
(note that any commutative algebra $C(X)$ is residually finite dimensional, independently on the topological dimension of $X$).
It is also clear that the existence of a spectral triple on $A$ satisfying \eqref{triv} implies that $A$ is residually finite
dimensional.

%\medskip \noindent \textbf{Acknowledgement.} The note was written during the visit of the authors at the Oberwolfach Institute
%(Research in Pairs programme).

\end{document}